\documentclass[11pt]{article}

\usepackage{amsmath}
\usepackage{amssymb}
\usepackage{eucal}

\begin{document}

\baselineskip 16pt

\title{On one generalization of finite nilpotent  groups}

\author{Zhang Chi \thanks{Research of the first author is supported by 
 China Scholarship Council and NNSF of 
China(11771409)}\\
{\small Department of Mathematics, University of Science and
Technology of China,}\\ {\small Hefei 230026, P. R. China}\\
{\small E-mail:
zcqxj32@mail.ustc.edu.cn}\\ \\
{ Alexander  N. Skiba}\\
{\small Department of Mathematics and Technologies of Programming,
  Francisk Skorina Gomel State University,}\\
{\small Gomel 246019, Belarus}\\
{\small E-mail: alexander.skiba49@gmail.com}}

\date{}
\maketitle

\date{}
\maketitle

\begin{abstract}  Let  $\sigma =\{\sigma_{i} | i\in I\}$ be  a
partition of the set $\Bbb{P}$ of all primes and  $G$   a finite group.
 A chief factor $H/K$ of $G$ is
said to be \emph{$\sigma$-central}  if the semidirect
 product $(H/K)\rtimes (G/C_{G}(H/K))$ is
a  $\sigma_{i}$-group for some $i=i(H/K)$.   $G$ is called  \emph{$\sigma$-nilpotent}
 if every  chief factor of $G$ is $\sigma$-central.   We say that $G$ is 
 \emph{semi-${\sigma}$-nilpotent} (respectively \emph{weakly semi-${\sigma}$-nilpotent}) if  
the normalizer $N_{G}(A)$ of every non-normal  (respectively  every non-subnormal)
 $\sigma$-nilpotent subgroup $A$  of $G$ is $\sigma$-nilpotent. 

In this paper we determine the structure of finite   
semi-${\sigma}$-nilpotent and weakly  semi-${\sigma}$-nilpotent groups.

\end{abstract}


\let\thefootnote\thefootnoteorig

\section{Introduction}

Throughout this paper, all groups are finite and $G$ always denotes
a finite group.  Moreover,
 $\mathbb{P}$ is the set of all    primes,
  $\pi \subseteq  \Bbb{P}$ and  $\pi' =  \Bbb{P} \setminus \pi$. If $n$ is an integer,
 the symbol $\pi (n)$
denotes the
 set of all primes dividing $n$; as usual,  $\pi (G)=\pi (|G|)$, the set of all
 primes dividing the order of $G$.

In what follows, $\sigma =\{\sigma_{i} | i\in I\}$ is  some
partition of $\Bbb{P}$, that is,
$\Bbb{P}=\bigcup_{i\in I} \sigma_{i}$ and $\sigma_{i}\cap
\sigma_{j}= \emptyset  $ for all $i\ne j$.  
By the analogy with the notation   $\pi (n)$, we write  $\sigma (n)$ to denote 
the set  $\{\sigma_{i} |\sigma_{i}\cap \pi (n)\ne 
 \emptyset  \}$;   $\sigma (G)=\sigma (|G|)$. A group
 is said to be  \emph{$\sigma$-primary}
 \cite{1} if   it  is a $\sigma _{i}$-group for some $i$.

  A chief factor $H/K$ of $G$ is
said to be \emph{$\sigma$-central} (in $G$) \cite{1} if the semidirect
 product $(H/K)\rtimes (G/C_{G}(H/K))$ is
$\sigma$-primary. The normal subgroup $E$ of $G$ is called  
\emph{$\sigma$-hypercentral} in $G$ if either $E=1$ or every chief factor of $G$ below $E$ is  
$\sigma$-central.
                        
Recall also that $G$ is called  \emph{$\sigma$-nilpotent} \cite{1} 
 if every  chief factor of $G$ is $\sigma$-central. 

An arbitrary group $G$ has two canonical $\sigma$-nilpotent subgroups of  particular 
importance in this context. The first of these is the 
\emph{$\sigma$-Fitting subgroup } $F_{\sigma}(G)$, that is, the product 
of all normal $\sigma$-nilpotent subgroups of $G$.  The other useful 
subgroup is the \emph{$\sigma$-hypercentre $Z_{\sigma}(G)$ of $G$}, that 
is, the  product of all $\sigma$-hypercentral subgroups of $G$. 

Note that in the classical case, 
when $\sigma = \sigma ^{1}=\{\{2\}, \{3\}, \ldots \}$ (we use here the notation in 
\cite{alg12}),   $F_{\sigma}(G)=F(G)$ is the Fitting subgroup 
 and $Z_{\sigma}(G)=Z_{\infty}(G)$ is the hypercentre of $G$.

In fact,  the $\sigma$-nilpotent groups are exactly the  groups $G$ which can be written in the
 form $G=G_{1} \times \cdots \times G_{t}$ for some $\sigma$-primary groups 
 $G_{1}, \ldots , G_{t}$ \cite{1}, and such groups have   proved to be very useful in
 the formation theory (see, in particular, the papers    \cite{19, 
20} and  the books \cite[Ch. IV]{bookShem}, \cite[Ch. 6]{15}). 
  In the recent years, the $\sigma$-nilpotent groups   have found new and to
 some extent unexpected  
applications in the theories of permutable and generalized subnormal 
subgroups  (see, in particular,  \cite{1, alg12},    \cite{3}--\cite{6} and the survey \cite{comm}).

In view of the results in the paper  \cite{belon}, the $\sigma$-nilpotent groups can
 be characterized  as the groups in which 
the normalizer of any $\sigma$-nilpotent subgroup is $\sigma$-nilpotent. 
Groups in which normalizers of all non-normal $\sigma$-nilpotent subgroups are
 $\sigma$-nilpotent may  be 
non-$\sigma$-nilpotent (see Example 1.3 below), and in the case when $\sigma = \sigma ^{1}$
 such  groups have been described in \cite[Ch. 4, Section 7]{We} (see also \cite{Sah}). 
In this paper, we determine the structure of  such groups $G$    for the case arbitrary
 $\sigma$.

{\bf Definition 1.1.}    We say that $G$ is 
(i) \emph{semi-${\sigma}$-nilpotent} if  
the normalizer of every non-normal $\sigma$-nilpotent subgroup 
of $G$ is $\sigma$-nilpotent;

(ii) \emph{weakly semi-${\sigma}$-nilpotent} if  
the normalizer of every non-subnormal $\sigma$-nilpotent subgroup 
of $G$ is $\sigma$-nilpotent;

(iii) \emph{weakly semi-nilpotent}  if $G$ is  weakly 
semi-${\sigma}^{1}$-nilpotent. 

{\bf Remark 1.2.}  (i)  Every ${\sigma}$-nilpotent group is  
semi-${\sigma}$-nilpotent, and every semi-${\sigma}$-nilpotent group is weakly 
semi-${\sigma}$-nilpotent.

(ii) The semi-${\sigma}^{1}$-nilpotent groups 
are exactly the \emph{semi-nilpotent groups} studied in   
\cite[Ch. 4, Section 7]{We}   (see also \cite{Sah}).

 (iii)   We show that $G$ is (weakly) semi-${\sigma}$-nilpotent if and only if   
the normalizer of every  non-normal (respectively non-subnormal) $\sigma$-primary  subgroup  
of $G$ is $\sigma$-nilpotent.  Since  every  $\sigma$-primary  group is
 $\sigma$-nilpotent, it is enough to show that 
if  the normalizer of every non-normal (respectively non-subnormal)  $\sigma$-primary
  subgroup  $A$ 
of $G$ is $\sigma$-nilpotent, then $G$ is ${\sigma}$-semi-nilpotent (respectively weakly
 semi-${\sigma}$-nilpotent). 
First note that   $A\ne 1$ and $A=A_{1} \times \cdots \times A_{n}$,
 where  $\{A_{1},  \ldots ,  A_{n}\}$  is a complete Hall $\sigma$-set of 
$A$. The subgroups  $A_{i}$ are  
characteristic in $A$, so  $N_{G}(A)=N_{G}(A_{1})
 \cap \cdots \cap N_{G}(A_{n})$, where either $N_{G}(A_{n})=G$ or $N_{G}(A_{n})$ is
 $\sigma$-nilpotent.  Since   $A$ is non-normal (respectively  non-subnormal) in $G$, 
there is $i$ such that $N_{G}(A_{n})$ is
 $\sigma$-nilpotent. Therefore $N_{G}(A)$ is
 $\sigma$-nilpotent by Lemma 2.2(i) below. Hence  $G$ is semi-${\sigma}$-nilpotent
 (respectively weakly  semi-${\sigma}$-nilpotent). 

{\bf Example 1.3.}    Let $p > q > r > t >  2$ be   primes, where $q$ divides $p-1$ and
 $t$ divides $r-1$,
  and let  $\sigma =\{\{p\},  \{q\},  \{p, q\}'\}$. Let 
 $R$ 
be  the quaternion group of order 8, $A$  a  group of order $p$, and let
 $B=C_{p}\rtimes C_{q}$ be
 a non-nilpotent  group  of order $pq$ and $C$ 
 a non-nilpotent  group  of order $rt$.   Then $B\times R$ is  a  non-$\sigma$-nilpotent
 semi-${\sigma}$-nilpotent group  and  $B \times C$ is not  semi-${\sigma}$-nilpotent.

Now let $G=A\times (Q\rtimes R)$, where $Q$ is a simple    ${\mathbb F}_{q}R$-module  which is    faithful  
for $R$.     Then for every subgroup $V$ of $R$
 we have $N_{G}(V)=A\times R$, so 
$G$ is  weakly  
semi-${\sigma}$-nilpotent. On the other hand, $QV$ is supersoluble for 
every  subgroup $V$ of $R$ of order 2  and so for some subgroup $L$ of $Q$ with 
$1 < L < Q$ we have  $V \leq N_{G}(L)$ and $[L, V]\ne 1$. Hence 
 $G$ is not semi-${\sigma}$-nilpotent.

Recall that  $G^{{\mathfrak{N}}_{\sigma}}$
  is the \emph{$\sigma$-nilpotent residual of $G$}, that is,  the 
intersection of all normal subgroups $N$ of $G$ with  $\sigma$-nilpotent 
quotient $G/N$.

 Our goal here is to determine the structure of  weakly  
semi-${\sigma}$-nilpotent and 
semi-${\sigma}$-nilpotent groups. In fact, the following concept is an 
important tool to achieve such a goal.

{\bf Definition 1.4.}     Let  $H$ be a
 ${\sigma}$-nilpotent subgroup of $G$.  Then we say that $H$ is 
  \emph{$\sigma$-Carter subgroup} of $G$ if $H$ is an \emph{${\mathfrak{N}}_{\sigma}$-covering 
subgroup of $G$} \cite[p. 101]{15}, that is, $U^{{\mathfrak{N}}_{\sigma}}H=U$
 for every subgroup $U$ of $G$  containing $H$. 

Note that in Example 1.3,  the subgroup $C_{q}C$
  is a $\sigma$-Carter subgroup of  the group $B \times C$.
It is clear also  that a group $H$ of a soluble group $G$ is a Carter subgroup of 
$G$ if and only if it is a $\sigma ^{1}$-Carter subgroup of $G$.

 A \emph{complete  set of Sylow subgroups of $G$}  
contains exactly one Sylow $p$-subgroup for each prime $p$ dividing $|G|$.
 In general, we say that a   set  ${\cal H}$ of subgroups of $G$ is a
 \emph{complete Hall $\sigma $-set} of $G$ \cite{2, comm}  if
 every member $\ne 1$ of  ${\cal H}$ is a Hall $\sigma _{i}$-subgroup of $G$
 for some $i$ and ${\cal H}$ contains exactly one Hall
 $\sigma _{i}$-subgroup of $G$ for every  $\sigma _{i}\in  \sigma (G)$.

Our first  result is the following

{\bf Theorem  A.} {\sl If  $G$ is  weakly semi-${\sigma}$-nilpotent, then:}

(i) {\sl $G$ has a complete Hall $\sigma$-set $\{H_{1}, \ldots , H_{t}\}$  
such that for some $1\leq r \leq t$ the subgroups 
$H_{1}, \ldots ,  H_{r}$ are normal in $G$,  $H_{i}$ is not   normal in $G$ for all $i > r$,
 and $$\langle H_{r+1},
 \ldots , H_{t} \rangle =H_{r+1}\times  \cdots \times H_{t}.$$}

(ii) {\sl If $G$ is not ${\sigma}$-nilpotent, then  $N_{G}(H_{i})$
 is a $\sigma$-Carter subgroup of $G$ for all $i > r$.}

(iii) {\sl $F_{\sigma}(G)$ is a maximal
 ${\sigma}$-nilpotent subgroup of $G$  and 
  $F_{\sigma}(G)=F_{0\sigma }(G)Z_{\sigma}(G)$, where
 $F_{0\sigma }(G)=H_{1} \cdots   H_{r}$.}

(iv) {\sl $V_{G}= Z_{\sigma}(G)$ 
for every maximal ${\sigma}$-nilpotent subgroup $V$ of $G$ such that  $G=F_{\sigma}(G)V$.
}

(v) {\sl   $G/F(G)$ is $\sigma$-nilpotent. }

On the basis of Theorem A we prove also the following

{\bf Theorem  B.}  {\sl Suppose that  $G$ is  semi-${\sigma}$-nilpotent,  
and let  $\{H_{1}, \ldots , H_{t}\}$ be a complete Hall $\sigma$-set of $G$, where  
 $H_{1}, \ldots ,  H_{r}$ are normal in $G$ and  $H_{i}$ is not   normal in 
$G$ for all $i > r$. Suppose also that  non-normal Sylow subgroups   of 
any  Schmidt subgroup $A\leq H_{i}$ 
have prime order for all $i > r$. 
Then:}

(i) {\sl   $G/F_{\sigma}(G)$ is abelian. }  

(ii) {\sl  If 
 $U$ is  any  maximal $\sigma$-nilpotent non-normal subgroup of $G$, then 
$U$ is a $\sigma$-Carter subgroup of $G$ and $U_{G}= Z_{\sigma}(G)$.}

(iii) {\sl If  the subgroups   $H_{1}, \ldots ,  H_{r}$ are nilpotent, then $G/F_{\sigma}(G)$
 is cyclic.  }

(iv) {\sl  Every quotient and every subgroup of $G$ are  semi-${\sigma}$-nilpotent}.
 
Now we consider some of corollaries of Theorems A and B in   the three classical cases.  
First of all note that  in the  case 
when $\sigma = \sigma ^{1}$, Theorems A and B not only cover the main results in  
\cite[Ch. 5 Section 7]{We} but they also  give the alternative proofs of  them.
  Moreover, in this case we get   from the theorems the following
 results.

{\bf  Corollary 1.4.} {\sl If  $G$ is  weakly semi-nilpotent, then:}
                 
(i) {\sl $G$ has a   complete  set of Sylow subgroups $\{P_{1}, \ldots , P_{t}\}$  such 
that for some $1\leq r \leq t$ the subgroups 
$P_{1}, \ldots ,  P_{r}$ are normal in $G$,  $P_{i}$ is not   normal in $G$ for all $i > r$,
 and $\langle P_{r+1},
 \ldots , P_{t} \rangle =P_{r+1}\times  \cdots \times P_{t}.$}

(ii) {\sl $F(G)$ is a maximal
 nilpotent subgroup of $G$  and 
  $F(G)=F_{0\sigma }(G)Z_{\infty}(G)$, where
 $F_{0\sigma }(G)=P_{1} \cdots   P_{r}$.}

(iii) {\sl If $G$ is not nilpotent, then $N_{G}(P_{i})$
 is a Carter subgroup of $G$ for all $ i > r$.}

{\bf Corollary 1.5} (See   Theorem 7.6 in \cite[Ch. 4]{We}).  {\sl If $G$ is 
 semi-nilpotent and $F_{0}(G)$ denotes the  product of its normal
 Sylow subgroups, then 
$G/F_{0}(G)$ is nilpotent. }

{\bf Corollary 1.6} (See   Theorem 7.8 in \cite[Ch. 4]{We}).  {\sl If $G$ is 
semi-nilpotent, then: }
 
(a) {\sl $F(G)$ is a maximal nilpotent subgroup of $G$.}

(b) {\sl If $U$ is  a maximal nilpotent subgroup of $G$ and $U$ is not normal in $G$, then
 $U_{G}=Z_{\infty}(G)$.}

{\bf Corollary 1.7} (See  Theorem 7.10 in \cite[Ch. 4]{We}). 
 {\sl The class of all semi-nilpotent groups is closed under taking subgroups
 and homomorphic images.}

 In the   other classical case when  $\sigma =\sigma ^{\pi}=\{\pi, 
\pi'\}$, 
 $G$  is  $\sigma ^{\pi}$-nilpotent 
 if and only if $G$ is  \emph{$\pi$-decomposable}, that is,
 $G=O_{\pi}(G)\times O_{\pi'}(G)$.

Thus $G$ is   semi-${\sigma}^{\pi}$-nilpotent if and only if the 
normalizer of every $\pi$-decomposable non-normal subgroup of $G$ is 
$\pi$-decomposable; $G$ is  weakly semi-${\sigma}^{\pi}$-nilpotent if and only if the 
normalizer of every $\pi$-decomposable non-subnormal subgroup of $G$ is 
$\pi$-decomposable.   
Therefore in this case we get from Theorems A and B the following  results.

{\bf Corollary 1.8.} {\sl Suppose that $G$ is not $\pi$-decomposable.  If   the 
normalizer of every $\pi$-decomposable non-subnormal subgroup of $G$ is 
$\pi$-decomposable, then:}

(i) {\sl $G$ has a  Hall  $\pi$-subgroup $H_{1}$   and a  Hall  
$\pi'$-subgroup $H_{2}$,  
 and exactly  one of these subgroups,  $H_{1}$ say,  is normal in $G$.}

(ii) {\sl   $G/F(G)$ is $\pi$-decomposable.  }

(iii)   {\sl $N_{G}(H_{2})$  is an
 $\mathfrak{F}$-covering  subgroup of $G$, where $\mathfrak{F}$ is the class
 of all $\pi$-decomposable groups.}

(iv) {\sl $O_{\pi}(G)\times O_{\pi'}(G)=H_{1}\times O_{\pi'}(G)$ is a maximal  
$\pi$-decomposable subgroup of $G$   and every element of $G$
 induces a $\pi'$-automorphism   on  every   chief factor  of $G$ 
  below  $O_{\pi'}(G)$.}

{\bf  Corollary 1.9.} {\sl Suppose that $G$ is not $\pi'$-closed and     the 
normalizer of every $\pi$-decomposable non-normal subgroup of $G$ is 
$\pi$-decomposable.  Then $G=H_{1}\rtimes H_{2}$, where  $H_{1}$ is a Hall 
$\pi$-subgroup  and    $H_{2}$  is a  Hall  
$\pi'$-subgroup of $G$. Moreover, if   non-normal Sylow subgroups   of 
any  Schmidt subgroup  $A\leq H_{2}$  
have prime order, then:}

(i) {\sl   $G/O_{\pi}(G)\times O_{\pi'}(G)$ is abelian. }

(ii) {\sl  Every  
   maximal $\pi$-decomposable non-normal subgroup of $G$ is an
 $\mathfrak{F}$-covering  subgroup of $G$, where $\mathfrak{F}$ is the class
 of all $\pi$-decomposable groups.}

(iii) {\sl  If $H_{1}$ is nilpotent, then 
$G/O_{\pi}(G)\times O_{\pi'}(G) $ is cyclic. }

  In fact, in the theory of   $\pi$-soluble groups
 ($\pi= \{p_{1}, \ldots , p_{n}\}$) we deal with the partition 
   $\sigma =\sigma ^{1\pi }=\{\{p_{1}\}, \ldots , \{p_{n}\}, \pi'\}$.  
Moreover,   $G$ is  $\sigma ^{1\pi }$-nilpotent
 if and only if $G$ is    \emph{$\pi$-special} \cite{Cun2},
 that is, $G=O_{p_{1}}(G)\times \cdots \times
 O_{p_{n}}(G)\times O_{\pi'}(G)$.
 
Thus $G$ is   semi-${\sigma}^{1\pi}$-nilpotent if and only if the 
normalizer of every $\pi$-special non-normal subgroup of $G$ is 
$\pi$-special; $G$ is  weakly semi-${\sigma}^{1\pi}$-nilpotent if and only if the 
normalizer of every $\pi$-special non-subnormal subgroup of $G$ is 
$\pi$-special.  
Therefore in this case we get from Theorems A and B the following  results.

 {\bf Corollary 1.10.} {\sl Let $P_{i}$ be a
 Sylow $p_{i}$-subgroup of $G$ for all $p\in \pi= \{p_{1}, \ldots , p_{n}\}$. 
If    the 
normalizer of every $\pi$-special non-subnormal subgroup of $G$ is 
$\pi$-special, then: }

(i) {\sl $G$ has a  Hall  $\pi'$-subgroup $H$
 and at least one of  subgroups  $P_{1}, \ldots , P_{n}, H$  is normal in $G$.}

(ii) {\sl $O_{p_{1}}(G)\times \cdots \times
 O_{p_{n}}(G)\times O_{\pi'}(G)$ is a maximal  
$\pi$-special subgroup of $G$.}

(iii) {\sl   $G/F(G)$ is $\pi$-special.   }

{\bf  Corollary 1.11.} {\sl Suppose that the  
normalizer of every $\pi$-special non-normal subgroup of $G$ is 
$\pi$-special. If    non-normal Sylow subgroups   of 
any  Schmidt $\pi'$-subgroup  of $G$  
have prime order, then:}

(i) {\sl   $G/(O_{p_{1}}(G)\times \cdots \times
 O_{p_{n}}(G)\times O_{\pi'}(G))$ is abelian. }

(ii) {\sl  Every  
   maximal $\pi$-special  non-normal subgroup of $G$ is an
 $\mathfrak{F}$-covering  subgroup of $G$, where $\mathfrak{F}$ is the class
 of all $\pi$-special groups.}

(iii) {\sl  If every normal in $G$ subgroup $A\in \{P_{1}, \ldots , P_{n}, H\}$ is nilpotent,
 then  $G/(O_{p_{1}}(G)\times \cdots \times
 O_{p_{n}}(G)\times O_{\pi'}(G))$ is cyclic. }

 \section{Preliminaries}

Recall that $G$ is said to be:  a \emph{$D_{\pi}$-group} if $G$ possesses a Hall 
$\pi$-subgroup $E$ and every  $\pi$-subgroup of $G$ is contained in some 
conjugate of $E$;  a \emph{$\sigma$-full group
 of Sylow type}  \cite{1} if every subgroup $E$ of $G$ is a $D_{\sigma _{i}}$-group for every
$\sigma _{i}\in \sigma (E)$; \emph{$\sigma$-soluble} \cite{1}  if
  every chief factor of $G$ is $\sigma$-primary.

{\bf Lemma 2.1 } (See Theorem A and B in \cite{2}).
   {\sl If $G$ is $\sigma$-soluble, then   $G$ is a 
$\sigma$-full group  of Sylow type and,   for every  $i$, $G$  has
  a Hall $\sigma _{i}'$-subgroup and every two
 Hall $\sigma _{i}'$-subgroups of $G$ are conjugate. }

A subgroup $A$ of $G$ is said to be \emph{${\sigma}$-subnormal}
  in $G$ \cite{1} 
 if   there is a subgroup chain  $A=A_{0} \leq A_{1} \leq \cdots \leq
A_{n}=G$  such that  either $A_{i-1}\trianglelefteq A_{i}$ or 
$A_{i}/(A_{i-1})_{A_{i}}$ is   $\sigma$-primary 
  for all $i=1, \ldots , n$.  Note that a subgroup $A$ of $G$ is subnormal 
in $G$ if and only if $A$ is ${\sigma}^{1}$-subnormal in $G$ (where
 ${\sigma}^{1}=\{\{2\}, \{3\}, \ldots \}$).

{\bf Lemma 2.2. }  (i)  {\sl 
 The class  of all $\sigma$-nilpotent groups ${\mathfrak{N}}_{\sigma}$ is closed under taking direct 
products, homomorphic images and  subgroups. Moreover, if  $H$ is a normal 
subgroup of $G$ and  $H/H\cap \Phi (G)$ is $\sigma$-nilpotent, then 
$H$ is $\sigma$-nilpotent   } (See Lemma 2.5 in \cite{2}).

(ii)   {\sl  $G$  is $\sigma $-nilpotent if and
 only if every subgroup of $G$ is ${\sigma}$-subnormal in $G$ } (See 
\cite[Proposition 3.4]{6}).

(iii)   {\sl  $G$  is $\sigma $-nilpotent if and  only if
 $G=G_{1}\times \cdots \times G_{n}$ for some $\sigma$-primary groups $G_{1}, \ldots ,
  G_{n}$ } (See 
\cite[Proposition 3.4]{6}).

{\bf Lemma 2.3} (See Lemma 2.6 in \cite{1}). {\sl Let  $A$,  $K$ and $N$ be subgroups of  $G$.
 Suppose that   $A$
is $\sigma$-subnormal in $G$ and $N$ is normal in $G$.  }

(1) {\sl If $N\leq K$ and $K/N$ is $\sigma$-subnormal in $G/N$, then $K$
is $\sigma$-subnormal in $G$}.

(2) {\sl $A\cap K$    is  $\sigma$-subnormal in   $K$}.
 
(3) {\sl If $A$  is $\sigma $-nilpotent, then  $A\leq F_{\sigma}(G)$.}

(4)  {\sl $AN/N$    is  $\sigma$-subnormal in   $G/N$}.

(5) {\sl If $A$  is a  Hall $\sigma _{i}$-subgroup of $G$ for some $i$, then  $A$ is normal in $G$.}

In view of Proposition 2.2.8  in \cite{15},
 we get from Lemma 2.2 the following

{\bf Lemma 2.4.} {\sl    If $N$ is a normal subgroup of $G$, then
 $(G/N)^{{\frak{N}}_{\sigma}}=G^{{\frak{N}}_{\sigma}}N/N.$  }

{\bf Lemma   2.5.} {\sl If  $G$ is  ${\sigma}$-soluble   and, for some 
$i$ and some   Hall $\sigma _{i}$-subgroup $H$ of $G$, $N_{G}(H)$
  is ${\sigma}$-nilpotent, then  $N_{G}(H)$ is a $\sigma$-Carter subgroup of $G$. }

{\bf Proof.}      Let $N=N_{G}(H)$ and $N\leq U\leq G$. 
Suppose that $U^{{\mathfrak{N}}_{\sigma}}N\ne U$ and let $M$ be a maximal 
subgroup of $U$ such that $U^{{\mathfrak{N}}_{\sigma}}N\leq M$. Then $M$ 
is $\sigma$-subnormal in $U$ by  Lemmas 2.2(i, ii) and 2.3(1), so    $U/M_{U}$ is a 
$\sigma _{j}$-group  for some $j$ since $U$ is  clearly ${\sigma}$-soluble.
Therefore   $|U:M|$ is a $\sigma _{j}$-number, 
so $j\ne i$ and hence $H\leq M_{U}$. But then $U=M_{U}N_{U}(H)\leq M < U$ by Lemma 2.1 and 
 the Frattini argument. This contradiction completes the proof of the 
lemma.

It is clear that if $A$ is $\sigma$-Carter subgroup of $G$, then $A$ is 
a 
$\sigma$-Carter subgroup in every subgroup of $G$ containing $A$. 
Moreover, in view of Proposition 2.3.14 in \cite{15}, the following useful facts are  true.

{\bf Lemma  2.6.} {\sl Let $H$ and  $R$ be subgroups of $G$, 
where $R$ is normal in $G$.  }

(i)  {\sl If
 $H$ is  a $\sigma$-Carter subgroup of $G$, then $HR/R$ is 
 a $\sigma$-Carter subgroup of $G/R$. }

(ii) {\sl If $U/R$ is 
 a $\sigma$-Carter subgroup of $G/R$ and $H$ is 
 a $\sigma$-Carter subgroup of $U$,  then 
 $H$ is  a $\sigma$-Carter subgroup of $G$. }

{\bf Lemma  2.7.} {\sl Suppose that   $G$ 
possesses a $\sigma$-Carter subgroup. If   $G$ is    ${\sigma}$-soluble,  then any two
 of its 
$\sigma$-Carter subgroups are conjugate. }

{\bf Proof.}  Assume that this lemma is false and let $G$ be a 
counterexample of minimal order. Then $|\sigma (G)| > 1$.

 Let   $A$ and
 $B$ be $\sigma$-Carter subgroups of $G$, and let 
 $R$  be a minimal normal subgroup of   $G$.  Then  $AR/R$  and
 $BR/R$ are  $\sigma$-Carter subgroups of $G/R$ by Lemma 2.6(i).
Therefore for some $x\in G$ we have $AR/R=B^{x}R/R$ by the choice of $G$. 
If $AR\ne G$, then $A$ and $B^{x}$ are conjugate in $AR$ by the choice of 
$G$ and so $A$ and $B$ are conjugate.

Now assume that    $AR=G=B^{x}R=BR$. If $R\leq A$, then $A=G$ is 
$\sigma$-nilpotent and so $A=B$. Therefore we can assume that  
$A_{G}=1=B_{G}$.

 Since $G$ is ${\sigma}$-soluble, $R$ is a ${\sigma}_{i}$-group for 
some $i$. Let $H$ be a Hall ${\sigma}_{i}'$-subgroup of $A$.
 Since $|\sigma (G)| > 1$, it follows that $H\ne 1$ and so  
$N=N_{G}(H)\ne 1$.   Since $A$ and
 $B$ be $\sigma$-Carter subgroups of $G$, both these subgroups 
are $\sigma$-nilpotent. Hence  $A\leq N$ and, for some $x\in G$,  $B^{x}\leq N$ by Lemma 2.1.
But then the choice of $G$ implies that  $A$ and $B^{x}$ are conjugate in $N$.
So we  again get that  $A$ and $B$ are conjugate.  The lemma is proved.

 If $G   \not \in {\mathfrak{N}}_{\sigma}$ but 
every proper  subgroup of $G$ belongs to  ${\mathfrak{N}}_{\sigma}$, then $G$ is called 
 an \emph{${\mathfrak{N}}_{\sigma}$-critical 
}  or a \emph{minimal non-$\sigma$-nilpotent} group. If $G$ is  an 
  ${\mathfrak{N}}_{{\sigma}^{1}}$-critical group, that is, $G$ is not nilpotent but
 every proper subgroup of $G$ is nilpotent, then   $G$ is said to be a \emph{Schmidt 
group}.

 {\bf Lemma 2.8} (See  \cite[Ch. V, Theorem 26.1]{bookShem}).  {\sl If  
$G$ is a   
 Schmidt group, then $G=P\rtimes Q$, where    $P=G^{\frak{N}}=G'$ 
 is a Sylow $p$-subgroup of $G$ and $Q=\langle x \rangle $ is a cyclic
 Sylow $q$-subgroup of $G$ with $\langle x^{q} \rangle \leq Z(G)\cap \Phi (G)$. Hence $Q^{G}=G$.    }

{\bf Lemma 2.9.} {\sl If 
   $G$ is an   ${\frak{N}}_{\sigma}$-critical 
group, then  $G$ is a Schmidt group.}

{\bf Proof.} For some $i$, $G$ is an    ${\frak{N}}_{\sigma _{0}}$-critical 
group, where  $\sigma _{0}=\{\sigma_{i}, \sigma_{i}'\}$. Hence  $G$ is a Schmidt group by 
\cite{belon}.

{\bf Lemma 2.10. }  {\sl  Let   $Z=Z_{\sigma}(G)$.
 Let  $A$, $B$   and $N$ be subgroups of $G$, where
$N$ is normal in $G$.}

(i) {\sl $Z$ is   ${\sigma}$-hypercentral in $G$. }

(ii) {\sl If $ N\leq Z$, then $Z/N= Z_{\sigma}(G/N)$.}

(iii) {\sl  $Z_{\sigma}(B)\cap A\leq Z_{\sigma}(B\cap A)$. }

(iv)   {\sl If $A$ is $\sigma$-nilpotent, then $ZA$   is also $\sigma$-nilpotent. Hence
$Z$ is contained in each maximal $\sigma$-nilpotent subgroup of $G$.}

(v) {\sl If $G/Z$ is $\sigma$-nilpotent, then $G$   is also $\sigma$-nilpotent.}

 {\bf Proof. } (i)   It is enough to consider the case when
$Z=A_{1}A_{2}$, where
$A_{1}$ and $A_{2}$ are normal ${\sigma}$-hypercentral subgroups of
 $G$. Moreover, in view of the  Jordan-H\"{o}lder theorem for the chief series,
it is enough to show that if   $A_{1}\leq K < H \leq A_{1}A_{2}$, then $H/K$
   is $\sigma$-central. But in this case we
have  $H=A_{1}(H\cap A_{2})$, where  $H\cap A_{2}\nleq  K$ and so 
from  
the $G$-isomorphism $(H\cap A_{2})/(K\cap A_{2})\simeq  (H\cap 
A_{2})K/K=H/K$ we get that $C_{G}(H/K)=C_{G}((H\cap A_{2})/(K\cap A_{2}))$
  and hence   $H/K$   is $\sigma$-central in $G$.

(ii) This  assertion is a  corollary  of Part (i) and
  the  Jordan-H\"{o}lder theorem for the chief series.

(iii)  First assume that $B=G$, and let  $1= Z_{0} < Z_{1}
< \cdots <  Z_{t} = Z$ be a chief
 series of $G$      below $Z$ and  $C_{i}=
C_{G}(Z_{i}/Z_{i-1})$. Now  consider the series
 $$1= Z_{0}\cap A \leq Z_{1}\cap A
 \leq \cdots  \leq  Z_{t} \cap A= Z\cap A.$$
We can assume without loss of generality that this series is a chief
series of $A$ below $Z\cap A$.

 Let   $i\in \{1,  \ldots , t \}$. Then, by Part (i),
$Z_{i}/Z_{i-1} $ is $\sigma$-central in $G$,
  $(Z_{i}/Z_{i-1})\rtimes (G/C_{i})$  is a $\sigma _{k}$-group say.
Hence    $(Z_{i}\cap A)/(Z_{i-1}\cap A)$ is a  $\sigma _{k}$-group. On the
other hand,
 $A/A\cap C_{i}\simeq C_{i}A/C_{i}$ is a  $\sigma _{k}$-group and   $$A\cap C_{i}\leq
 C_{A}((Z_{i}\cap A)/(Z_{i-1}\cap A)).$$ Thus  $(Z_{i}\cap
A)/(Z_{i-1}\cap A)$
is $\sigma$-central in $A$.   Therefore, in view of the  Jordan-H\"{o}lder theorem
for the chief series, we have  $Z\cap A\leq Z_{\sigma}(A)$.

Now  assume that  $B$ is  any subgroup of $G$. Then,
 in view of  the preceding paragraph,   we have
 $$ Z_{\sigma}(B) \cap A = Z_{\sigma}(B) \cap
(B\cap A)\leq Z_{\sigma}(B\cap A).$$

(iv)   Since $A$ is $\sigma$-nilpotent, $ZA/Z\simeq A/A\cap Z$ is
$\sigma$-nilpotent by Lemma 2.2(i). On the other hand, $Z\leq Z_{\sigma}(ZA)$
 by Part (iii).   Hence $ZA$ is $\sigma$-nilpotent  by Part (i).

(v) This assertion follows from Part (i).

The lemma  is proved.

 The following lemma is a corollary of Lemmas 2.2(i) and 2.10(v).

 {\bf Lemma 2.11.} {\sl  $F_{\sigma}(G)/\Phi (G)=F_{\sigma}(G/\Phi(G))$ and  
$F_{\sigma}(G)/Z_{\sigma}(G)=F_{\sigma}(G/Z_{\sigma}(G))$.}

\section{Proofs of the main results}

{\bf Proof of Theorem A.}  Assume that this theorem is false and let $G$ be a 
counterexample of minimal order.  Then $G$ is not $\sigma$-nilpotent.

(1) {\sl   Every proper subgroup  $E$ of $G$ is  weakly 
semi-${\sigma}$-nilpotent. Hence the conclusion of the theorem   holds for $E$.}

  Let $V$ be a non-subnormal   $\sigma$-nilpotent subgroup
 of $E$. 
Then $V$ is not subnormal   
 in $
G$ by Lemma 2.3(2),  so $N_{G}(V)$ is 
 $\sigma$-nilpotent by 
hypothesis. Hence  $N_{E}(V)=N_{G}(V)\cap E$ is  $\sigma$-nilpotent by Lemma 2.2(i).

(2) {\sl  Every proper quotient $G/N$ of $G$ (that is, $N\ne 1$)
  is  weakly semi-${\sigma}$-nilpotent. Hence the conclusion of the theorem   holds for  
   $G/N$.   }

In view of Remark 1.2(iii)   and the choice of $G$, it is enough to show that if $U/N$ is any
 non-subnormal  
  $\sigma$-primary 
subgroup of $G/N$, then  $N_{G/N}(U/N)$ is $\sigma$-nilpotent. 
 We can 
assume   without loss of generality  that $N$ is a minimal normal subgroup of $G$.

Since  $U/N$ is not 
 subnormal    
  in  $G/N$, $U/N  < G/N$ and $U$  is not 
 subnormal  in $G$. Hence $U$ is a proper 
subgroup of $G$, which implies that $U$ is $\sigma$-soluble  by Claim (1).
Hence   $N$ is a $\sigma _{i}$-group for some $i$.

 If $U/N$ is 
a $\sigma _{i}$-group, then  $U$ is $\sigma$-primary and so $N_{G}(U)$ is  
 $\sigma$-nilpotent.
 Hence $N_{G/N}(U/N)=N_{G}(U)/N$ is $\sigma$-nilpotent
 by Lemma 2.2(i).  Now suppose that 
$U/N$ is 
a $\sigma _{j}$-group for some $j\ne i$. Then $N$ has a  complement $V$ in 
$U$ by the Schur-Zassenhaus theorem. Moreover,  from the Feit-Thompson  
theorem it follows that at least one of the groups $N$ or $U/N$ is soluble 
and so every two complements to $N$ in $U$ are conjugate in $U$. Therefore 
$N_{G}(U)=N_{G}(NV)=NN_{G}(V)$. Since $U=NV$ is not 
 subnormal  in $G$,     $V$   is not 
 subnormal in $G$ by Lemma 2.3(1, 4)  
  and  so $N_{G}(V)$ is $\sigma$-nilpotent. Hence
  $N_{G/N}(U/N)=N_{G}(U)/N$ is $\sigma$-nilpotent.

(3) {\sl If    $A$ is an ${\mathfrak{N}}_{\sigma}$-critical 
subgroup of $G$, then  $A=P\rtimes Q$, where    $P=A^{\frak{N}}=A'$ 
 is a Sylow $p$-subgroup of $A$ and $Q$ is a 
 Sylow $q$-subgroup of $A$ for some different primes $p$ and $q$.  
Moreover,  $P$ is subnormal in $G$ and so $P\leq O_{p}(G)$.  }

The first assertion of the claim directly follows from 
 Lemmas 2.8 and 2.9. Since $A$ is not 
$\sigma$-nilpotent, $P$ is subnormal in $G$ by hypothesis. Therefore 
$ P\leq O_{p}(G)$ by Lemma 2.3(3).

(4)  {\sl $G$  is $\sigma$-soluble.}

  Suppose that this is false. Then  $G$ is a non-abelian simple group since every
 proper section of $G$ is $\sigma$-soluble by Claims (1) and (2).   Moreover, 
 $G$ is not $\sigma$-nilpotent and so  it has  an ${\mathfrak{N}}_{\sigma}$-critical 
subgroup  $A$.   Claim (3) implies that for some Sylow subgroup $P$ of $A$ 
we have   $1 <  P\leq O_{p}(G) < G$. 
    This contradiction shows that  we have (4).

 (5)  {\sl Statements (i) and (ii) hold for $G$.}

 Since $G$ is  $\sigma$-soluble by Claim (4), 
it is a $\sigma$-full group of Sylow type by Lemma 2.1. In particular, $G$ 
possesses a complete   Hall $\sigma$-set $\{H_{1}, \ldots , H_{t}\}$. Then 
there is an index $k$ such that $H_{k}$ is not subnormal  in $G$ 
  by Lemma 2.3(5) 
 since $G$ is not $\sigma$-nilpotent. Then    $N_{G}(H_{k})$
 is $\sigma$-nilpotent by hypothesis, so 
  $N_{G}(H_{i})$ is  a $\sigma$-Carter subgroup of $G$ by Lemma 2.5
 for all $i  > r$.

If for some $j\ne k$ the subgroup $H_{j}$ is   not subnormal in $G$, then 
$N_{G}(H_{j})$ is also a $\sigma$-Carter subgroup of $G$. But then  
$N_{G}(H_{k})$  and $N_{G}(H_{j})$  are conjugate in $G$ by Lemma 2.7. 
Hence for some $x\in G$ we have  $H_{k}^{x}\leq  N_{G}(H_{j})$. Therefore, 
since $G$ is not $\sigma$-nilpotent, 
there is a complete Hall $\sigma$-set  $\{L_{1}, \ldots , L_{t}\}$ of $G$ 
such that for some 
 $1\leq r < t$ the subgroups $L _{1}, \ldots , L_{r}$ are normal in $G$, 
 $L_{i}$ is not   normal in $G$ for all $i > r$,  
 and
 $\langle L_{r+1},
 \ldots , L_{t} \rangle =L_{r+1}\times  \cdots \times L_{t}$.

(6) {\sl Every subgroup $V$
 of $G$ containing $ F_{\sigma}(G)$ is $\sigma$-subnormal in $G$, so  
 $F_{\sigma}(V)=F_{\sigma}(G)$.}

From Claim (5) it follows that $H_{1},   \ldots ,  H_{r} \leq F_{\sigma}(G)$ and 
$$G/F_{\sigma}(G)=F_{\sigma}(G)(H_{r+1}\times  \cdots \times 
H_{t})/F_{\sigma}(G)\simeq (H_{r+1}\times  \cdots \times 
H_{t})/((H_{r+1}\times  \cdots \times 
H_{t}) \cap  F_{\sigma}(G))$$ is  ${\sigma}$-nilpotent.  Hence every 
subgroup of $G/F_{\sigma}(G)$ is ${\sigma}$-subnormal in   
$G/F_{\sigma}(G)$ by Lemma 2.2(ii).  Therefore $V$ is  
${\sigma}$-subnormal in  $G$ by Lemma 2.3(1), so
 $F_{\sigma}(V)\leq F_{\sigma}(G)\leq F_{\sigma}(V)$ by Lemma 2.3(3). Hence 
we have (6).

(7)  {\sl Statement (iii)  holds for $G$.}

First note that $F_{\sigma}(G)$ is a maximal
 ${\sigma}$-nilpotent subgroup of $G$ by Claim (6).    
In fact, $F_{\sigma}(G)=F_{0\sigma }(G)\times O_{\sigma _{i_{1}}}(G)\times \cdots \times
 O_{\sigma _{i_{m}}}(G)$ for some $i_{1}, \ldots , i_{m} \subseteq \{r+1, 
\ldots , t\}$.
Moreover, in view of Claim (5), we get  clearly that $G/C_{G}(O_{\sigma _{i_{k}}}(G))$  is a 
$\sigma _{i_{k}}$-group and so    $O_{\sigma _{i_{k}}}(G)\leq Z_{\sigma}(G)$. 
  Hence 
  $F_{\sigma}(G)=F_{0\sigma }(G)Z_{\sigma}(G)$.

(8)  {\sl Statement (iv) holds for $G$.}

First we show that $U_{G}\leq  Z_{\sigma}(G)$ for every ${\sigma}$-nilpotent subgroup $U$
 of $G$               such that  $G=F_{\sigma}(G)U$.  Suppose that this is false. Then $U_{G}\ne 1$. Let $R$ be a minimal 
normal subgroup of $G$ contained in $U$ and $C=C_{G}(R)$. 
  Then $$G/R=(F_{\sigma}(G)R/R)(U/R)=F_{\sigma}(G/R)(U/R),$$ so $$U_{G}/R=(U/R)_{G/R}\leq 
 Z_{\sigma}(G/R)$$ by Claim (2). Since $G$ is $\sigma$-soluble, $R$ 
is a $\sigma _{i}$-group for some $i$.  Moreover, from $G=F_{\sigma}(G)U$ 
and Lemma 2.1 we get that for some Hall $\sigma _{i}'$-subgroups $E$, $V$ and 
$W$ of $G$, of $F_{\sigma}(G)$ and of  $U$, respectively, we have $E=VW$. 
But $R\leq F_{\sigma}(G)\cap U$, where $F_{\sigma}(G)$ and $U$ are  
$\sigma$-nilpotent.   Therefore $E\leq C$, so $R/1$ is $\sigma$-central in 
$G$. Hence $R\leq Z_{\sigma}(G)$ and so $Z_{\sigma}(G/R)=Z_{\sigma}(G)/R$ by Lemma 2.10(ii). 
But then $U_{G}\leq Z_{\sigma}(G)$. Finally, if $V$ is any   
maximal ${\sigma}$-nilpotent subgroup  of $G$ with $G=F_{\sigma}(G)V$, then  $Z_{\sigma}(G)\leq  V$ by 
Lemma 2.11(iv) and so $V_{G}= Z_{\sigma}(G)$.

(9)   {\sl Statement (v) holds for $G$.}

In view of Lemma 2.2(i), it is enough to show that $D=G^{{\mathfrak{N}}_{\sigma}}$ is 
nilpotent. Assume that this is false. Then $D\ne 1$, and for   
any  minimal normal subgroup $R$  of $G$ we have  that 
 $(G/R)^{{\mathfrak{N}}_{\sigma}}=RD/R \simeq D/D\cap R $ is nilpotent by 
Claim (2) and  Lemmas 
2.2(i) and 2.4. Moreover, Lemma 2.2(i) implies 
that $R$ is a unique minimal normal subgroup of $G$, $R\leq D$ and $R\nleq 
\Phi (G)$.    
Since  $G$ is not ${\sigma}$-nilpotent, Claim (3) and  \cite[Ch. A, 15.6]{DH} imply that 
$R=C_{G}(R)=O_{p}(G)=F(G)$  
 for some prime $p$. 
Then  $R  < D$ and $G=R\rtimes M$, where $M$ is not $\sigma$-nilpotent, 
and so  
$M$ has an  ${\mathfrak{N}}_{\sigma}$-critical 
subgroup   $A$.   Claim (3) implies  that    for some prime $q$ dividing $|A|$ and for a Sylow 
$q$-subgroup $Q$ of $A$ we have $1   <  Q\leq F(G)\cap M=R\cap M=1$. 
   This contradiction completes the proof of (9).

From Claims (5), (7), (8) and (9) it follows that  the conclusion of the theorem is 
true for $G$, contrary to the choice of $G$. The theorem is proved.

{\bf Proof of Theorem B.}
  Assume that this theorem is false and let $G$ be a 
counterexample of minimal order. Then $G$ is not ${\sigma}$-nilpotent. Nevertheless, $G$  is 
${\sigma}$-soluble by Theorem A. 
Let  $F_{0\sigma }(G)=H_{1} \cdots   H_{r}$ and 
$ E=H_{r+1} \cdots H_{t}$.   Then $E$ is $\sigma$-nilpotent  by Theorem A(ii).

(1) {\sl   Every proper subgroup  $E$ of $G$ is   
semi-${\sigma}$-nilpotent.  
 Hence  Statements (i) and (ii) hold for $E$} (See Claim (1) in the proof of Theorem A).

(2) {\sl   The hypothesis holds for every  proper quotient $G/N$ of $G$. 
    Hence  Statements (i),  (ii) and (iv) hold for 
$G/N$. }

It is not difficult to show  that  $G/N$  is   semi-${\sigma}$-nilpotent (see
 Claim (2) in the proof of Theorem A).

Now let  $U/N$ be any Schmidt $\sigma _{i}$-subgroup  of $G/N$  such that 
$U/N\leq W/N$ for some non-normal in $G/N$   Hall $\sigma _{i}$-subgroup $W/N$ of $G/N$.  
In view of Lemma 2.1, we can assume without loss of generality that 
$W/N=H_{i}N/N$.
 Let $L$ be
 any minimal supplement to $N$ in $U$.   Then $L\cap N\leq \Phi (L)$ and,  by Lemma 2.8,    
 $U/N=LN/N\simeq L/L\cap N$  is a $\sigma _{i}$-group and $L/L\cap N
  =(P/L\cap N)\rtimes (Q/L\cap N)$, where 
   $P/L\cap N=(L/L\cap N)^{\frak{N}}=(L/L\cap N)'$ 
 is a Sylow $p$-subgroup of $L/L\cap N$ and $Q/L\cap N=\langle x \rangle $ is a cyclic
 Sylow $q$-subgroup of $L/L\cap N$ with
 $V/L\cap N=\langle x^{q} \rangle =\Phi (Q/L\cap N)\leq \Phi
 (L/L\cap N)\cap  Z(L/L\cap N)$ 
 and $p, q \in \sigma _{i}$.  Suppose that 
$|Q/L\cap N|  > q$.   Then  $L\cap N <  V$.

In view of Lemma 2.2(i), a Sylow $p$-subgroup of $L$ is normal in $L$. 
 Hence, in view of   Lemma 
2.8, for any Schmidt subgroup $A$ of $L$ we have 
$A =A_{p}\rtimes A_{q}$, where  $A_{p}$ is a Sylow $p$-subgroup of $A$,   
 $ A_{q}$  is a Sylow $q$-subgroup of $A$ and $(A_{q})^{A}=A$. 
  We can assume without loss of 
generality that  $A_{q}(L\cap N)/(L\cap N) \leq Q/L\cap N$. 
Then    $A_{q}(L\cap N)/(L\cap N)
 \nleq V/L\cap N$ since $V\leq \Phi (L)$.
It follows that $A_{q}\nleq N$.  Since $W/N=H_{i}N/N$ is not normal in 
$G/N$, $H_{i}$ is not normal in $G$. But for some $x\in G$ we have $A^{x}\leq 
H_{i}$, so  $|A_{q}^{x}|=|A_{q}|=q$ by hypothesis.

Note that    $|Q/V|=q$ since $Q/L\cap N$
 is cyclic and $V/L\cap N=\Phi (Q/L\cap N)$. Hence  
 $$(V/L\cap N)(A_{q}(L\cap N)/(L\cap N))=(V/L\cap N)\times
 (A_{q}(L\cap N)/(L\cap N))=Q/(L\cap N) ,$$
 which implies that $Q/(L\cap N)$ is not cyclic. This contradiction shows  that $|Q/L\cap N|=q$, so 
for a Sylow $q$-subgroup 
$S$ of $U/N$ we have $|S|=q$.  Therefore the hypothesis holds for  $G/N$.  
Hence we have (2) by the choice of $G$

(3) {\sl If    $A$ is an ${\mathfrak{N}}_{\sigma}$-critical 
subgroup of $G$, then  $A=P\rtimes Q$, where    $P=A^{\frak{N}}=A'$ 
 is a Sylow $p$-subgroup of $A$ and $Q$ is a 
 Sylow $q$-subgroup of $A$ for some different primes $p$ and $q$. Moreover, the
 subgroup $P$ is normal in $G$.  Hence $G$  has  an  abelian minimal normal subgroup $R$}
(See Claim (3) in the proof of Theorem A).

(4)  {\sl Statement (i) holds for $G$.}

  In view of Lemma 
2.2(i), it is enough to show that $G'$ is  ${\sigma}$-nilpotent. Suppose 
that this is false.

(a) {\sl $R=C_{G}(R)=O_{p}(G)=F(G)\nleq \Phi (G)$ for some prime $p$ and $|R|  > p$}.

 From Claim (2) it follows that for every minimal normal subgroup $N$ of $G$,
 $(G/N)'=G'N/N\simeq G'/G'\cap N$
 is $\sigma$-nilpotent.  If $R\ne N$, it follows that  $G'/((G'\cap N)\cap 
(G'\cap R))=G'/1$ is  $\sigma$-nilpotent by Lemma 2.2(i). Therefore $R$ is a unique 
minimal normal subgroup of $G$, $R\leq D$ and $R\nleq \Phi (G)$ by Lemma 
2.2(i). Hence 
$R=C_{G}(R)=O_{p}(G)=F(G)$   by    Theorem 15.6 in \cite[Ch. A]{DH}, so
  $|R| > p$ since otherwise   
       $G/R=G/C_{G}(R)$ is cyclic, which implies that   $G'=R$ is  ${\sigma}$-nilpotent.

(b) {\sl 
 $F_{\sigma}(V)=F_{\sigma}(G)$ for every subgroup $V$
 of $G$ containing $ F_{\sigma}(G)$}  (See Claim (6) in the proof of Theorem A).

(c)  {\sl $G=H_{1} \rtimes H_{2}$, where $R\leq H_{1}=F_{\sigma}(G)$ and $H_{2}$ is
 a minimal non-abelian group.}

From    Theorem A and  Claim  (a) it follows that $r=1$  and $R\leq H_{1}=F_{\sigma}(G)$.

 Now let    
$W=F_{\sigma}(G)V$, where $V$ is a  maximal subgroup of $E$.  Then 
 $F_{\sigma}(G)=F_{\sigma}(W)$ by Claim (b), so  $W/F_{\sigma}(W)=W/F_{\sigma}(G)\simeq V$ is 
abelian by Claim (1). Therefore $E$  is not abelian but   
 every proper subgroup of $E$  is abelian, so  $E=H_{2}$ since $E$ is
 $\sigma$-nilpotent.  Hence we have (c).

(d)  {\sl $H_{1}=R$  is a Sylow $p$-subgroup of $G$ and every  subgroup $H\ne 1$ of
 $H_{2}$ acts irreducibly  on $R$. Hence every proper  subgroup $H$  of $H_{2}$ is cyclic.}

 Suppose that $|\pi (H_{1})| > 1$.  There is a Sylow $p$-subgroup $P$ 
of $H_{1}$ such that $H_{2}\leq N_{G}(P)$ by Claim (c) and the Frattini argument.
Let $K=PH_{2}$. Then $K  < G$ and  $P=H_{1}\cap K$ is normal in $K$, so
 $R\leq P=F_{\sigma}(K)$ since  $C_{G}(R)=R$ by Claim (a). 
Then $K/F_{\sigma}(K)=K/P\simeq H_{2}$ is abelian by Claim (1), a contradiction.
Hence $H_{1}$  is a normal Sylow $p$-subgroup of $G$. Hence 
$H_{1}\leq F(G)\leq C_{G}(R)=R$ by \cite[Ch. A, 13.8(b)]{DH}, so $H_{1}=R$.

Now let $S=RH$.  By 
the Maschke theorem, $R=R_{1}\times \cdots \times R_{n}$, where $R_{i}$ 
is a minimal normal subgroup of $S$ for all $i$. Then $R=C_{S}(R)=C_{S}(R_{1}) \cap 
\cdots \cap C_{S}(R_{n})$. Hence, for some $i$,
 the subgroup $R_{i}H$ is not $\sigma$-nilpotent and so it has 
 an ${\mathfrak{N}}_{\sigma}$-critical 
subgroup  $A$ such that $1  < A'$ is normal in $G$ by Claim (3). But then
 $R\leq A$. Therefore $i=1$, so we have (d) since $H$ is abelian by Claim (c).

(e) {\sl $H_{2}$ is not nilpotent. Hence $ |\pi (H_{2})| > 1$.}

Suppose that $H_{2}=Q\times H$ is nilpotent, where $Q\ne 1$ is a Sylow $q$-subgroup of $H_{2}$.
If $H\ne 1$, then $Q$ and $H$ are proper subgroups of $H_{2}$ and so the groups $Q$,
 $H$ and $H_{2}$ are abelian by Claim (c). Therefore  $H_{2}=Q$ is a $q$-group. Then, 
since every maximal subgroup of  $H_{2}$ is cyclic by Claim (d), $q=2$ by
 \cite[Ch. 5, Theorems 4.3, 4.4]{Gor}. Therefore $|R|=p$, contrary to 
Claim (a).  Hence we have (e).

(f)  {\sl  $H_{2}=A\rtimes B$, where $A=C_{H_{2}}(A)$ is
 a group of prime order $q\ne p$ and $B=\langle a \rangle$ is 
 a  group of order $r$ for some prime $r\not \in \{p, q\}$. } 

From Claims 
 (d) and (e)  it follows that $H_{2}$ is a Schmidt group with cyclic
 Sylow subgroups. Therefore  Claim (f) follows from the hypothesis and  Lemma 2.8.

 {\sl   Final contradiction for (4). }      Suppose that for some $x=yz\in RA$,
 where $y\in R$ and $z\in A$, 
 we have $xa=ax$. Then  $x\in N_{G}(B)$, so $R\cap \langle x \rangle =1$ 
 since $B$ acts irreducible on $R$ by Claim (d). 
Hence $\langle x \rangle $ is a $q$-group and $V=\langle x \rangle B$ is abelian group such that 
$B\cap R=1$. Hence from the isomorphism $G/R\simeq H_{2}$ we get that $x=1$.
  Therefore  $a$ induces a fixed-point-free  automorphism    on $RA$ and hence $RA$
 is nilpotent by the
 Thompson theorem \cite[Ch. 10, Theorem 2.1]{Gor}. But then $A\leq 
C_{G}(R)=R$. This contradiction completes the proof of (4).

(5)  {\sl Statement (ii) holds for $G$.}

 Suppose that this is false. 
By Lemma 2.10(iv),
  $Z_{\sigma}(G)\leq  U$. On the other,  
$U/Z_{\sigma}(G)$ is a
maximal $\sigma$-nilpotent non-normal subgroup of $G/Z_{\sigma}(G)$  by 
Lemma 2.10(v). Hence 
in the case $Z_{\sigma}(G)\ne 1$ Claim (2) implies that  
$U/Z_{\sigma}(G)$ is a  $\sigma$-Carter subgroup  $G/Z_{\sigma}(G)$,  so 
 $U$ is a  $\sigma$-Carter subgroup  of $G$ by Lemma 2.6(ii).  Hence   $Z_{\sigma}(G)=1$, so   
    Theorem A(iii) implies that $F_{\sigma}(G)= 
F_{0\sigma}(G)=H_{1}\cdots H_{r}$. Hence $E\simeq  G/F_{0\sigma}(G)$ is abelian by Claim (4).

Let $V=F_{\sigma}(G)U $. 
If $V=G$, then for some $x$ we have $H_{r+1}^{x}\leq U$ by Lemma 2.1. Hence   
$U\leq N_{G}(H_{r+1}^{x})$ and so $U= N_{G}(H_{r+1}^{x})$  is a 
 $\sigma$-Carter subgroup of $G$ by   Theorem A(ii).  
Therefore  
$V=F_{\sigma}(G)U $ is a normal proper subgroup of  $G$. 
Let  $x\in G$. If the  subgroup $U^{x}$ is normal in $V$, then  $U^{x}$ is 
subnormal in $G$  and so $U^{x}, U\leq F_{\sigma}(G) $ by Lemma 2.3(3), which implies that
 $U=F_{\sigma}(G) $ is normal in $G$ since $F_{\sigma}(G) $ and $U$ are
 maximal $\sigma$-nilpotent subgroups of $G$ by Theorem A(iii). This contradiction shows 
that  $U^{x}$ and  $U$ are non-normal maximal $\sigma$-nilpotent subgroups 
of  $V$. Since $  V  < G$, Claim (1) implies that  $U^{x}$ and  $U$ are 
$\sigma$-Carter  subgroups  of  $V$.  
 Since $V$ is $\sigma$-soluble, $U$ and $U^{x}$ are conjugate in $V$ 
by Lemma 2.7. Therefore $G=VN_{G}(U)$ by the Frattini argument.  Since $U$ 
is a maximal 
$\sigma$-nilpotent non-normal subgroup  of $G$, $U=N_{G}(U)$. Hence 
$G=VU=(F_{\sigma}(G)U)U =F_{\sigma}(G)U < G$. This contradiction completes
 the proof of the fact that  every  maximal $\sigma$-nilpotent non-normal subgroup $U$
  of $G$ is a $\sigma$-Carter subgroup of $G$. But then   
$G=F_{\sigma}(G)U$ since $G/F_{\sigma}(G)$ is  $\sigma$-nilpotent by Claim (4) and so  
   $U_{G}= Z_{\sigma}(G)$ by Theorem A(iv).  Hence we have (5).

(6) {\sl If $F_{0\sigma} (G)\leq F(G)$, then $G/F_{\sigma}(G)$ is cyclic.}

Assume that this  is false.

(i) {\sl  $\Phi (F_{0\sigma}(G))=1$. Hence $F_{0\sigma} (G)$ is the direct product of some 
minimal normal subgroups $R_{1}, \ldots , R_{k}$ of $G$.}

Suppose that  $\Phi (F_{0\sigma} (G))\ne 1$ and let $N$ be a minimal normal
 subgroup of $G$ contained in 
$\Phi (F_{0\sigma} (G))\leq \Phi (G)$. Then $N$ is a $p$-group for some prime $p$. 

We show that the hypothesis holds for $G/N$. First note that $G/N$ is 
semi-${\sigma}$-nilpotent  by Claim (2). Now let $V/N$ be a normal  Hall
 $\sigma _{i}$-subgroup of $G/N$ for some $\sigma _{i}\in \sigma (G/N)$. 
If $p\in \sigma _{i}$, then $V$ is normal Hall $\sigma _{i}$-subgroup of $G$, so
 $V\leq F(G)$ by hypothesis 
and hence $V/N\leq F(G)N/N\leq F(G/N)$. Now assume that $p\not \in \sigma _{i}$
 and let $W$ be  a 
  Hall  $\sigma _{i}$-subgroup of $V$. Then $W$ is   a 
  Hall  $\sigma _{i}$-subgroup of $G$.  Moreover, every two   
  Hall  $\sigma _{i}$-subgroups of $V$ are conjugate in $V$ by Lemma 2.1, so
 $G=VN_{G}(W)=NWN_{G}(W)=NN_{G}(W)=N_{G}(W)$ by the Frattini argument.
 Therefore $W\leq F(G)$, so $V/N=WN/N\leq F(G/N)$. 
Hence $F_{0\sigma}(G/N)\leq F(G/N)$, so the hypothesis holds for $G/N$.  
The choice of $G$ and Lemma 2.11 imply that  $(G/N)/F_{\sigma}(G/N) =
 (G/N)/(F_{\sigma}(G)/N)\simeq G/F_{\sigma}(G)$
 is cyclic, a contradiction.  Hence $\Phi (F_{0\sigma}(G))=1$, so we have (i) by 
\cite[Ch. A, Theorem 10.6(c)]{DH}.

(ii)   {\sl $Z_{\sigma}(G)=1$.  Hence $F_{0\sigma} (G)=F_{\sigma} (G)=F(G)$.}

Since $Z_{\sigma}(G/Z_{\sigma}(G))=1$  by Lemma 2.10(ii), Lemma 2.11 and   Theorem A(iii)
 imply that  
 $$ 
F_{0\sigma}(G/Z_{\sigma}(G))=F_{\sigma}(G/Z_{\sigma}(G))
=F_{\sigma}(G)/Z_{\sigma}(G)=F_{0\sigma} 
(G)Z_{\sigma}(G)/Z_{\sigma}(G), $$ where
 $F_{0\sigma} (G)\leq F(G)$ and  so    
$F_{0\sigma}(G/Z_{\sigma}(G))\leq F(G/Z_{\sigma}(G))$.  Therefore the 
hypothesis holds for $G/Z_{\sigma}(G)$ and hence, in the case when $Z_{\sigma}(G)\ne 1$, 
 $G/F_{\sigma}(G)\simeq 
(G/Z_{\sigma}(G))/F_{\sigma}(G/Z_{\sigma}(G))$ is cyclic by the choice of $G$.
Hence we have (ii).   

{\sl Final contradiction for (6).}    
 Since $E\simeq G/F(G)$ is abelian by Claims  (4) and  (ii) and   $G$ is not nilpotent,  there
 is an index $i$ such that $V=R_{i}\rtimes E$ is not 
nilpotent. Then $C_{R_{i}}(E)\ne R_{i}$.   By the Maschke theorem, $R_{i}= 
L_{1}\times \cdots \times L_{m}$
  for some minimal normal subgroups $L_{1}, \ldots ,  L_{m}$ of $V$. Then, since 
$C_{R_{i}}(E)\ne R_{i}$,  for some $j$ we have $L_{j}\rtimes E\ne L_{j}\times E$. Hence 
$L_{j}E$ contains a Schmidt subgroup $A_{p}\rtimes A_{q}$ such that $A_{p}=R_{i}$, so $m=1$. 
But then  $E$ acts irreducible 
on $R_{i}$ and  hence $G/F(G)\simeq E$ is cyclic.  This 
contradiction completes the proof of (6).

From Claims (1), (2), (4),  (5) and  (6) it follows that  the conclusion of the theorem is 
true for $G$, contrary to the choice of $G$. The theorem is proved.

 \end{document}